\input amstex
\input amsppt.sty
\hsize 30pc
\vsize 47pc
\def\nmb#1#2{#2}         
\def\cit#1#2{\ifx#1!\cite{#2}\else#2\fi} 
\def\totoc{}             
\def\idx{}               
\def\ign#1{}             
\redefine\o{\circ}

\define\al{\alpha}
\define\be{\beta}
\define\ga{\gamma}
\define\de{\delta}
\define\ep{\varepsilon}

\define\ph{\varphi}

\define\ps{\psi}

\define\Ga{\Gamma}
\define\De{\Delta}

\define\La{\Lambda}

\define\Ph{\Phi}

\redefine\i{^{-1}}
\define\x{\times}
\define\p{\partial}
\define\Id{\operatorname{Id}}
\def\today{\ifcase\month\or
 January\or February\or March\or April\or May\or June\or
 July\or August\or September\or October\or November\or December\fi
 \space\number\day, \number\year}
\topmatter
\title  Construction of completely integrable systems by
        Poisson mappings
\endtitle
\author
J\. Grabowski,
G\. Marmo,
P\. W\. Michor
\endauthor
\leftheadtext{\smc Grabowski, Marmo, Michor}
\rightheadtext{\smc Construction of completely integrable systems}
\affil
Erwin Schr\"odinger International Institute of Mathematical Physics,
Wien, Austria
\endaffil
\address J\. Grabowski: Institute of Mathematics,
University of Warsaw,
ul. Banacha 2, PL 02-097 Warsaw, Poland; and
Mathematical Institute, Polish Academy of Sciences,
ul. \'Sniadeckich 8, P.O. Box 137, PL 00-950 Warsaw, Poland
\endaddress
\email jagrab\@mimuw.edu.pl \endemail
\address
G. Marmo:
Dipart\. di Scienze Fisiche - Universit\`a di Napoli,
Mostra d'Oltremare, Pad.19, I-80125 Napoli, Italy.
\endaddress
\email gimarmo\@na.infn.it \endemail
\address
P\. W\. Michor: Institut f\"ur Mathematik, Universit\"at Wien,
Strudlhofgasse 4, A-1090 Wien, Austria; and
Erwin Schr\"odinger International Institute of Mathematical Physics,
Boltzmanngasse 9, A-1090 Wien, Austria
\endaddress
\email peter.michor\@univie.ac.at, peter.michor\@esi.ac.at \endemail
\date Sept. 01, 1999 
\enddate
\thanks J.~Grabowski was supported by KBN, grant Nr 2 P03A 031 17.  
\endthanks
\abstract Pulling back sets of functions in involution by Poisson
mappings and adding Casimir functions during the process allows 
to construct completely integrable systems. 
Some examples are investigated in detail.
\endabstract
\endtopmatter
\document

\head\totoc\nmb0{1}. Introduction \endhead

The standard notion of complete integrability is the so called
Liouville-Arnold integrability: a Hamiltonian system on a
$2n$-dimensional symplectic manifold $M$ is said to be completely
integrable if it has $n$ first integrals in involution which are
functionally independent on some open and dense subset of $M$.

It is natural to extend the notion of complete integrability to
sytems defined on Poisson manifolds $(N,\La)$ by requiring that on
each symplectic leaf such system defines a completely integrable system in
the usual sense. This  generalization implies that an integrable
system is associated to a maximal abelian Poisson subalgebra of
$(C^\infty(N),\{\quad,\quad\}_\La)$. The dynamical system $\Ga$
associated to a 1-form $\al$ on $N$ via $\Ga=i_\al\La$ will define a
Hamiltonian system on a symplectic leaf $S$ with embedding
$\ep_S:S\to N$ if we have $\ep_s^*\al=dH_S$.

If this is the case for any symplectic leaf we may write
$\al=\sum_kf_k\,dg^k$ where the $f_k$ are Casimir functions for $\La$.
When all the $g^k$'s belong to a sufficiently large set of functions
in involution which are functionally independent on each leaf, the
dynamical system is completely integrable. For some cases, even if
$\ep_S^*d\al\ne 0$, we get an integrable system in the generalized
sense of \cit!{1}.

Of course, integrable systems are not easy to find. Recently, in the
paper \cit!{3} we came accross a beautiful idea to construct
completely integrable systems by using coproducts in Poisson-Hopf
algebras. In this paper we put this construction into a geometric
perspective in order to understand better which are the essential
ideas that make the construction possible. In addition to this, we
construct a full family of Poisson-Hopf algebras associated with a
parametrized family of Poisson-Lie structures on the group
$SB(2,\Bbb C)$. The standard Lie-Poisson structures on $SB(2,\Bbb C)$
with $SU(2)$ and $SL(2,\Bbb R)$ as dual groups are included in this
scheme. This Poisson-Hopf algebras can be viewed as geometrical
version of the corresponding quantum groups --- deformations of the
universal enveloping algebra of the associated Lie algebras. We also
present symplectic realizations of the corresponding commutation
rules in the deformed algebras.

\head\totoc\nmb0{2}. Constructing integrable systems by Poisson maps
\endhead

\subhead\nmb.{2.1}. Complete integrability on Poisson manifolds
\endsubhead
If $(M,\La)$ is a Poisson manifold, a Hamiltonian system
$H\in C^\infty(M)$ is called \idx{\it completely integrable} if it
admits a complete set of first integrals in involution:
There are $f_1,\dots,f_k\in C^\infty(M)$ with $\{f_i,f_k\}_\La=0$ and
$\{f_i,H\}_\La=0$ such that on each symplectic
leaf (or an open dense set of
symplectic leaves) $H$ together with a suitable subset of
$f_1,\dots,f_k$ is Liouville-Arnold integrable.

\subhead\nmb.{2.2}. Constructing families of functions in involution
by Poisson maps \endsubhead
\newline
Let $\Ph_i: (M_{i+1},\La_{i+1})\to (M_i,\La_i)$ be Poisson maps
between manifolds, so that
$\{f,g\}_{i}\o\Ph_i=\{f\o\Ph_i,g\o\Ph_i\}_{i+1}$.
If we have a family of functions $\Cal F_1\subset C^\infty(M_1)$
in involution on $(M_1,\La_1)$, we may consider the family
$\Cal F_2=(\Cal F_1\o\Ph_1)\cup C_2\subset C^\infty(M_2)$ where $C_2$
is a complete set of Casimir functions on $(M_2,\La_2)$, and so on:
$$\CD
(M_1,\La_1) @. \hskip 2cm \Cal F_1, \text{ in involution}\\
@A{\Ph_1}AA     @. \\
(M_2,\La_2) @. \hskip 2cm \Cal F_2 = (\Cal F_1\o\Ph_1)\cup C_2 \\
@A{\Ph_2}AA     @. \\
(M_3,\La_3) @. \hskip 2cm \Cal F_3 = (\Cal F_2\o\Ph_2)\cup C_3 \\
@A{\Ph_3}AA     @. \\
..         @.
\endCD$$

\subhead\nmb.{2.3}. Poisson actions and multiplications \endsubhead
We shall apply the procedure of \nmb!{2.2} mainly in the following
situation: Consider $(M_1\x M_2,\La_1\x \La_2)$. Then for the
algebras of smooth functions we have
$C^\infty(M_1\x M_2)\cong C^\infty(M_1)\tilde\otimes C^\infty(M_2)$
for some suitable completed tensor product, where
$(f_1\otimes f_2)(x,y)=f_1(x)f_2(y)$. Then
$$
\{f_1\otimes f_2, g_1\otimes g_2\}_{\La_1\x \La_2}
     = \{f_1,g_1\}_{\La_1}\otimes g_1 g_2 + f_1 f_2
     \otimes \{g_1,g_2\}_{\La_1}.
$$
So if $c_1$ is a Casimir function of $(M_1,\La_1)$, then
$c_1\otimes 1$ is a Casimir function of $(M_1\x M_2,\La_1\x \La_2)$.
In this sense The Casimir functions of $(M_1,\La_1)$ and those of
$(M_2,\La_2)$ extend both to Casimir functions on
$(M_1\x M_2,\La_1\x \La_2)$.

If $\Ph:(M\x M,\La\x\La)\to (M,\La)$ is a Poisson map (for example
the multiplication of a Lie Poisson group) we may use it for the
procedure of \nmb!{2.2}. If $\Ph$ is associative then $\De_\Ph:
f\mapsto f\o \Ph$ is coassociative. But this is not essential for
applying the procedure in which $M_n=\prod^nM$ and $\Ph_n$ is a 
Cartesian product of $\Ph$ with identities. For example,
$$
M @<{\Ph}<< M\x M @<{\Ph\x \Id_M}<< M\x M\x M
     @<{\Id_M\x\Ph\x\Id_M}<< M\x M\x M\x M \gets \dots
$$
We start with a set of functions $\Cal F_1\subset C^\infty(M)$ in
involution and with a basis $\Cal C$ of all Casimirs.
Then $\Cal F_n \subset C^\infty(\prod^n M)$ is given
recursively by
$$
\Cal F_{n+1} =
   (\Cal F_{n}\o\Ph_n)\cup \{\Cal C\otimes 1\otimes\dots\otimes 1,
    1\otimes \Cal C\otimes 1\otimes \dots \otimes 1 , \dots\}
$$
and furnishes a family of functions in involution on $\prod^{n+1}M$.
If $\Ph$ is associative (so $\De_\Ph$ is coassociative) then the
result does not depend on the `path' chosen to define the $\Ph_n$'s.

Another possibility is to consider a Poisson mapping
$\Ph:(M\x N,\La_M\x\La_N)\to (N,\La_N)$ (for example
a Lie-Poisson action on $N$ of a Lie Poisson group $M$)
and to apply the procedure as follows:
$$
N @<{\Ph}<< M\x N @<{\Id_M\x\Ph}<< M\x M\x N
     @<{\Id_{M\x M}\x\Ph}<< M\x M\x M\x N \gets \dots
$$

\subhead\nmb.{2.4} \endsubhead
We may extend the procedure described in \nmb!{2.3} as follows.
We assume that we have furthermore Poisson manifolds (e.g\. symplectic ones)
$N_1,\dots N_n$ and Poisson mappings $\ph_i: N_i\to M$. The product
map
$\ph=\ph_1\x \dots\x \ph_n: N_1\x\dots \x N_n \to \prod^n M$ is a
Poisson map. Let $\Cal F_n$ be the set of functions in involution on
$\prod^n M$ constructed in \nmb!{2.3}. Then $\Cal F_n\o \ph$ is a set
of functions in involution on $\prod_i N_i$.

Standard examples of Poisson maps $\ph_i:N_i\to M$ are the canonical
embeddings of symplectic leaves $N_i$ of the Poisson manifold $M$. In
this case, the Casimir functions
$1\otimes 1\otimes \dots \otimes c\otimes \dots \otimes 1$
are constants on $N_1\x\dots \x N_n$, but the coproducts
$\De_\Ph c$, $(\De_\Ph\x\Id_M)\circ\De_\Ph c$, etc.,
are usually no longer Casimirs and hence sometimes
give rise to completely integrable systems on $N_1\x\dots\x N_n$.
See example \nmb!{3.1}.

\head\totoc\nmb0{3}. Examples \endhead

\subhead\nmb.{3.1}. Example \endsubhead
Let $M=\frak{su}(2)^*$ be the dual space of the Lie algebra
$\frak{su}(2)$. It carries a Kostant-Kirillov-Souriau Poisson
structure which is given in linear coordinates by
$$
\La= z\p_x\wedge \p_y + x\p_y\wedge \p_z + y\p_z\wedge \p_x.
$$
Since $\La$ is linear, we have the obvious Poisson map
$$
\Ph: M\x M\to M,\qquad \Ph(x_1,y_1,z_1,x_2,y_2,z_2) =
(x_1+x_2,y_1+y_2,z_1+z_2).
$$
A Casimir function for $\La$ is $c=x^2+y^2+z^2$. According to our
procedure in \nmb!{2.3} the functions
$$\align
c\o\Ph &= (x_1+x_2)^2 + (y_1+y_2)^2 + (z_1+z_2)^2, \\
c\otimes 1 &= x_1^2+y_1^2+z_1^2, \\
1\otimes c &= x_2^2+y_2^2+z_2^2, \\
f\o \Ph &= f(x_1+x_2,y_1+y_2,z_1+z_2),
\endalign$$
are functions in involution on $M\x M$,
where $f$ is an arbitrary function on $M$.
If we take for $N$ the symplectic leaf $N=c\i(1)$ which is a
2-dimensional sphere $S^2$, the Casimir functions $c\otimes 1$ and
$1\otimes c$ pull back to constants on $N\x N=S^2\x S^2$. However,
$c\o \Ph$ and $f\o \Ph$ are in involution and hence
$H=\frac12 (c\o \Ph) - 1 = x_1x_2+y_1y_2+z_1z_2 $ defines a
completely integrable system on the symplectic manifold
$N\x N\subset M\x M$. The system defined by the Hamiltonian function
$H$ is, in fact, completely integrable on each symplectic leaf of
$M\x M$, so that we get a completely integrable system on $M\x M$
whose dynamics is given by the vector field
$$\align
\Ga = (&z_2y_1-y_2z_1)\p_{x_1} + (z_1y_2-y_1z_2)\p_{x_2} +\\
      (&x_2z_1-z_2x_1)\p_{y_1} + (x_1z_2-z_1x_2)\p_{y_2} +\\
      (&y_2x_1-x_2y_1)\p_{z_1} + (y_1x_2-x_1y_2)\p_{z_2}.
\endalign$$
This vector field is tangent to all products of spheres
since
$c\otimes 1 = x_1^2+y_1^2+z_1^2$ and $1\otimes c = x_2^2+y_2^2+z_2^2$
are first integrals, and on $N\x N$ it induces the motion which can
be interpreted as associated with a `spin-spin'-interaction:
$$
\dot{\vec J_1} = \vec J_1 \times \vec J_2,\quad
\dot{\vec J_2} = \vec J_2 \times \vec J_1;
$$
The points on the spheres move in such a way that the velocity of
each of them is the vector product of the two position vectors.
Stationary solutions occupy the same or opposite points on the
sphere.

In sperical coordinates, the same system can be given a different
interpretation:
$$
z_i = \sin \be_i,\quad
y_i = \cos\be_i\;\sin\al_i,\quad
x_i = \cos\be_i\;\cos\al_i.
$$
In canonical coordinates $p_i=\sin\be_i, q_i = \al_i$ we get the
Hamiltonian function in the form
$$
H=p_1p_2 + \sqrt{(1-p_1^2)(1-p_2^2)}\cos(q_1-q_2).
$$
Since $H_1=\De(z^2)-\De(c)$ is in involution with $\De(z^2)$ we can
consider the completely integrable system given by the Hamiltonian
$$
H_1=p_1^2+p_2^2 -2 \sqrt{(1-p_1^2)(1-p_2^2)}\cos(q_1-q_2).
$$
Let us remark that our Hamiltonian $H$ is a slight modification of
the Hamiltonian
$$
H_0 = p_1p_2 - p_1p_2\cos(q_1-q_2)
$$
obtained in \cit!{3}.

We can inductively apply our procedure to get a completely integrable
system on $\prod^k M$ with Hamiltonian
$$
H^{(k)} = \sum_{i<j}^k (x_ix_j+y_iy_j+z_iz_j)
$$
which reduces in canonical coordinates on $\prod^k N$  to
$$
H^{(k)} = \sum_{i<j}^k \left(p_ip_j +
\sqrt{(1-p_i^2)(1-p_j^2)}\cos(q_i-q_j)\right).
$$

\subhead\nmb.{3.2}. Example \endsubhead
We consider the following symplectic realization \cit!{8} of the
Lie algebra $\frak{su}(2)$ in $T^*\Bbb R^2$:
$$
x = \tfrac12(q_1q_2+p_1p_2), \quad
y = \tfrac12(p_1q_2-q_1p_2), \quad
z = \tfrac14(p_1^2+q_1^2-p_2^2-q_2^2).
$$
This defines a Poisson morphism
$$
\ps: T^*\Bbb R^2 \to \frak{su}(2)^*,\quad
(q_1,q_2,p_1,p_2)\mapsto (x,y,z),
$$
which is the momentum map of the corresponding Hamiltonian action of
the group $SU(2)$.

As before, we consider the Casimir function
$c=x^2+y^2+z^2$ on $\frak{su}(2)$.
This time, however,
$$
F=c\o\ps=\tfrac1{16}\left(p_1^2+p_2^2+q_1^2+q_2^2\right)^2
$$
is not a Casimir function on the symplectic manifold $T^*\Bbb R^2$.
The functions in involution on
$\frak{su}(2)^*\x \frak{su}(2)^*$
from example \nmb!{3.1} give rise to functions in involution on
$T^*\Bbb R^2\x T^*\Bbb R^2= T^*\Bbb R^4$ as in \nmb!{2.4}, where
$\tilde q_1$ etc\. denote the functions on the second copy of
$T^*\Bbb R^2$:
$$
F_1 = F(q_1,q_2,p_1,p_2), \quad
F_2 = F(\tilde q_1,\tilde q_2,\tilde p_1,\tilde p_2), \quad
(\De c)\o \ps = F_1+F_2+H,
$$
where
$$\multline
H = \tfrac14\Bigl(
  (q_1q_2+p_1p_2)(\tilde q_1\tilde q_2+\tilde p_1\tilde p_2)
 +(p_1q_2-q_1p_2)(\tilde p_1\tilde q_2-\tilde q_1\tilde p_2)\Bigr)\\
 +\tfrac1{16}(p_1^2+q_1^2-p_2^2-q_2^2)
  (\tilde p_1^2+\tilde q_1^2-\tilde p_2^2-\tilde q_2^2)
\endmultline$$
and $(\De f)\o\ps = f(G_1,G_2,G_3)$, where
$$\align
G_1 &= \tfrac12(q_1q_2+p_1p_2 +
     \tilde q_1\tilde q_2+\tilde p_1\tilde p_2),\\
G_2 &= \tfrac12(p_1q_2-q_1p_2
     +\tilde p_1\tilde q_2-\tilde q_1\tilde p_2),\\
G_3 &= \tfrac14(p_1^2+q_1^2-p_2^2-q_2^2 +
     \tilde p_1^2+\tilde q_1^2-\tilde p_2^2-\tilde q_2^2).
\endalign$$
Hence, we have 4 independent functions in involution on $T^*\Bbb R^4$
which define completely integrable systems. As Hamiltonian functions
we can take the pure interaction term $H$. The trajectories of the
corresponding dynamics $\Ga_H$ lie on the intersections of the level
sets of $F_1$ and $F_2$ (which are, topologically, products of
3-dimensional spheres) and the level sets of $H$, and, say, $G_1$
(which are, generically, 4-dimensional tori).
Note that $G_2$ and $G_3$ are additional constants of the motion.
The whole set $\{F_1,F_2, H, G_1,G_2,G_3\}$ is, however, not
independent, since $G_1^2+G_2^2+G_3^2=F_1+F_2+H$.
The dynamics $\Ga_H$ on $T^*\Bbb R^4\cong\Bbb R^8$ is described by a
rather complicated vector field whose coefficients are polynomials of
degree 3.

The functions $G_1,G_2,G_3$ define the diagonal action of $SU(2)$ on
$T^*\Bbb R^2\x T^*\Bbb R^2$ which preserves $\Ga_H$. The dynamics on
$S^2\x S^2$ from example \nmb!{3.1} can be obtained via symplectic
reduction with respect to this  action.

\subhead\nmb.{3.3}. Example \endsubhead
Let us now consider the Lie group $M=SB(2,\Bbb C)$ of all matrices of
the form
$$
A=\pmatrix e^{-kz/2} & x+iy \\
           0         & e^{kz/2} \endpmatrix,
$$
where $k\ne0$ is fixed and $z,x,y\in \Bbb R$ may be viewed as global
coordinates on $SB(2,\Bbb C)$.
The coproduct corresponding to the group multiplication
$\Ph:SB(2,\Bbb C)\x SB(2,\Bbb C)\to SB(2,\Bbb C)$ is given by
$$\align
\De(z) &= z\otimes 1 + 1\otimes z, \\
\De(x) &= x\otimes e^{kz/2} + e^{-kz/2}\otimes x, \\
\De(y) &= y\otimes e^{kz/2} + e^{-kz/2}\otimes y.
\endalign$$
$\Ph$ is a Poisson map for each of the following Poisson structures
(parameterized by $\al,\be,\ga\in \Bbb R$):
$$
\{z,x\} = \be y, \quad
\{y,z\} = \al x, \quad
\{x,y\} = \ga\frac{\sin(kz)}{k}. \tag{*}
$$
For $\al=\be=\ga=1$ we get the `classsical realization' of the
quantum $SU(2)$ group, and for $\al=\be=-\ga=1$ we get the `classsical
realization' of the quantum $SL(2,\Bbb R)$ of Drinfeld and Jimbo,
\cit!{6}.

For $k\to 0$ we get the Lie algebras
$$
\{z,x\} = \be y, \quad
\{y,z\} = \al x, \quad
\{x,y\} = \ga z,
$$
with the standard cobrackets $\De(u)=u\otimes 1 + 1\otimes u$,
corresponding to the addition in $\frak g^*$.

A Casimir function for the Poisson bracket \thetag{*} is
$c=\al x^2+\be y^2 + \frac{4\ga}{k^2}\sinh^2(kz/2)$, which for
$k\to 0$ goes to $c_0=\al x^2 + \be y^2 + \ga z^2$, see \cit!{7}.

The generic symplectic leaves of the Poisson structure  are
2-dimensional (except for the trivial case $\al=\be=\ga=0$). As in
example \nmb!{3.1} the Hamiltonian $H=\frac12\De(c)$ defines a
completely integrable system on $SB(2,\Bbb C)\x SB(2,\Bbb C)$. In the
coordinates $x,y,z$ the Hamiltonian $H$ has the form
$$\multline
H= \tfrac12\left(\al x_1^2+\be y_1^2 +
     \tfrac{4\ga}{k^2}\sinh^2(\tfrac{kz_1}2)\right)e^{kz_2} \\
+ \tfrac12\left(\al x_2^2+\be y_2^2 +
     \tfrac{4\ga}{k^2}\sinh^2(\tfrac{kz_2}2)\right)e^{-kz_1} \\
+ \left(\al x_1x_2 + \be y_1y_2 +
     \tfrac{4\ga}{k^2}\sinh(\tfrac{kz_1}2)
     \sinh(\tfrac{kz_2}2)\right)e^{k(z_1-z_2)/2}.
\endmultline$$
In the limit for $k\to 0$ we get
$$
H_0 = c_0\otimes 1 + 1\otimes c_0 + (\al x_1x_2 + \be y_1y_2 + \ga z_1z_2)
$$
and for $\al=\be=\ga=1$ we are in the situation of example
\nmb!{3.1}.

In all generality, however, it is difficult to express the dynamics
explicitly since we deal simultaneously with a parametrized family of
structures for which even the topology of the symplectic leaves
changes.

Even in the case $\al=\be=\ga=1$, where it is known
\cit!{2} that $(SB(2,\Bbb C),\La)$ is equivalent, as a Poisson
manifold, to $\frak{su}(2)^*$ with the Kostant-Kirillov-Souriau
structure $\La_0$ described in example \nmb!{3.1}, the dynamics
described by $H$ in the deformed case may differ from that of example
\nmb!{3.1}. The reason is that $(SB(2,\Bbb C),\La)$ is not equivalent to
$(\frak{su}(2),\La_0)$ as a Lie Poisson group, since $SU(2,\Bbb C)$ is
not commutative. In particular, the deformed coproduct is not
cocommutative and the interaction we obtain is not symmetric.

In order to work in canonical coordinates let us introduce a
symplectic realization of the commutation rules \thetag * with
$\al=\de^2>0$ and $\be=1$:
$$
X = \sqrt{a^2-\tfrac{4\ga}{k^2}\sinh^2(\tfrac{kp}2)}\,
     \frac{\sin(\de q)}{\de}, \quad
Y= \sqrt{a^2-\tfrac{4\ga}{k^2}\sinh^2(\tfrac{kp}2)}\,
     \cos(\de q), \quad
Z= p,
$$
where $a\ge 0$ and $a>0$ if $\ga=0$.
In particular, if $\ga>0$ we get the deformed $SU(2)$, and if $\ga<0$
we get the deformed $SL(2,\Bbb R)$.
In this realization the Casimir function is $c=a^2$ and the
Hamiltonian reads
$$\multline
H = e^{k(p_2-p_1)/2}\biggl(\sqrt{
     \left(a^2-\tfrac{4\ga}{k^2}\sinh^2(\tfrac{kp_1}2)\right)
     \left(a^2-\tfrac{4\ga}{k^2}\sinh^2(\tfrac{kp_2}2)\right)
     }\,\cos(\de(q_1-q_2)) + \\
+ a^2\cosh(k\tfrac{p_1+p_2}{2})
     + \tfrac{4\ga}{k^2}\sinh(\tfrac{kp_1}2)\sinh(\tfrac{kp_2}2)
\biggr).
\endmultline$$
This Hamiltonian is quite complicated. But if we put $a=0$, $\ga=-1$,
and $\de=1$ we get
$$
H_1 = 4e^{k(p_2-p_1)/2} \tfrac{1}{k^2}\sinh(\tfrac{kp_1}2)
     \sinh(\tfrac{kp_2}2)(\cos(q_1-q_2) -1)
$$
which is the Hamiltonian obtained in \cit!{3} for the deformed
$SL(2,\Bbb R)$.

\subhead\nmb.{3.4}. Example  \endsubhead
A slight modification of the previous example which, at least
formally, is not dealing with Lie-Poisson groups, is the the
following. Let $M$ be the space of all upper triangular matrices of
the form
$$
A=\pmatrix a & x+iy \\
           0 & b \endpmatrix,
$$
where $a,b,x,y\in \Bbb R$. We use these as coordinates on $M$.
We consider the Poisson structure $\La$ on $M$ with Poisson bracket
$$\alignat3
\{x,a\} &= \be ya, &\quad \{x,b\} &=-\be yb, &\quad
     \{x,y\} &=-\ga (b^2-a^2), \tag{**}\\
\{y,a\} &= \al xa, &\quad \{y,b\} &= \al xb, &\quad
     \{a,b\} &= 0,
\endalignat$$
where $\al,\be,\ga\in\Bbb R$ paramerterize a family of Poisson
brackets. Matrix multiplication on $M$ leads to the coproduct:
$$
\De(a) = a\otimes a, \quad
\De(b) = b\otimes b, \quad
\De(x) = a\otimes x + x\otimes b, \quad
\De(y) = a\otimes y + y \otimes b.
$$
The matrix multiplication turns out to be a Poisson mapping with
respect to all brackets \thetag{**}. But $M$ is not a Lie-Poisson
group since it contains elements which are not invertible. On the
other hand, all $\La$ are tangent to $SB(2,\Bbb C)\subset M$ and give
there the brackets $\thetag{*}$ with slightly modified coeficients
$\al,\be,\ga$, if we parameterize $a=e^{-kz/2}$, $b=e^{kz/2}$.

The Poisson tensor has 2 independent Casimirs:
$c_1=ab$ and $c_2= \al x^2+\be y^2+\ga(a^2+b^2)$.
The symplectic leaves are, generically, 2-dimensional (we assume that
$\al^2+\be^2+\ga^2\ne0$), and as before,
$$\align
H=\tfrac12\De(c_2)
= {}&\tfrac12 a_1^2(\al x_2^2+\be y_2^2+\ga(a_2^2+b_2^2))+\\
&\tfrac12 b_2^2(\al x_1^2+\be y_1^2+\ga(a_1^2+b_1^2)) +\\
&a_1b_2(\al x_1x_2+\be y_1y_2-\ga a_1b_2)
\endalign$$
describes a completely integrable system on $M\x M$ which, on
$SB(2,\Bbb C)= c_1\i(1)$, coincides with a system from example
\nmb!{3.3}.

\subhead\nmb.{3.5}. Example \endsubhead
This is of different type.
Let $D=G.G^*=G^*.G$ be a complete Drinfeld double group. The
decompositions give us two Poisson projections
$\pi_1,\pi_2:D\to G^*$ onto the Lie-Poisson group $G^*$.
Let $\Cal F, \Cal F'$ be two families of functions in involution on
$G^*$. It is known that pull backs by $\pi_1$ and by $\pi_2$ commute
with respect to the symplectic structure on $D$.
So we can take $\Cal F_2 = (\Cal F_1\o\pi_1) \cup (\Cal F_2\o\pi_2)$
as a set of functions in involution on $D$.

For example,
$$
SL(2,\Bbb C) = SU(2).SB(2,\Bbb C) = SB(2,\Bbb C).SU(2),
$$
where $\pi_1$ denotes the projection onto the left factor
$SB(2,\Bbb C)$, and $\pi_2$ onto the the right one.
Let $\Cal F$ consist of the Casimir $C$ and some function $f$.
Then the functions $c\o\pi_1, f\o\pi_1, g\o\pi_2$ are in involution
on $D$, where $f$ and $g$ are arbitrary in $C^\infty(SB(2,\Bbb C))$,
and they are generically independent, so we have:
\block{\it
$H=f\o \pi_1$ is a completely integrable system on the symplectic
$SL(2,\Bbb C)$  for any $f\in C^\infty(SB(2,\Bbb C))$.
}\endblock

\subhead \nmb.{3.6}. Concluding remarks \endsubhead

1. We have shown that by using Poisson-compatible coproducts it is possible
to generate interacting systems while preserving the complete integrability.
We have given some examples. These can be extended to arbitrary Lie Poisson
pairs. 

2. The interaction we get is a 2-body interaction, one may wonder if it would
not be possible to obtain non-factorizable n-body interactions by using
n-ary brackets or n-ary operations (\cit!{7}, \cit!{8}, \cit!{9}).

3. The composition procedure does not use the existence of an inverse for
each element in the product, therefore the procedure may be extended to
any algebra. Is it possible to obtain interacting systems with fermionic
degrees of freedom by using graded algebras?

4. The procedure may clearly be extended to infinite dimensions. Can 
one use it to obtain interacting fields?

We shall come back to some of these questions in the near future.

\enddocument